\documentclass[12pt]{amsart}
\usepackage{amssymb}
\usepackage{graphicx}

\theoremstyle{plain}
\newtheorem{theorem}{Theorem}[section]
\newtheorem{lemma}[theorem]{Lemma}

\theoremstyle{definition}
\newtheorem{definition}[theorem]{Definition}
\theoremstyle{remark}

\newtheorem{remark}[theorem]{Remark}
\numberwithin{equation}{section}

  \renewcommand\epsilon\varepsilon
  \renewcommand\phi\varphi

\begin{document}
\author{Norayr Matevosyan and Peter Markowich}
\thanks{This work was supported by the RICAM (Austrian Academy od Sciences)
and by the WITTGENSTEIN AWARD 2000 of Peter Markowich, funded by the
Austrian Science Fund (FWF)}
\address{Peter Markowich\\
Institute of Mathematics \\
University of Vienna\\
Boltzmanngasse 9\\
A-1090 Vienna }
\email{peter.markowich@univie.ac.at}
\address{Norayr Matevosyan\\
Johann Radon Institute for Computational and Applied Mathematics (RICAM)\\
Austrian Academy of Sciences\\
Altenbergerstraße 69 \\
A-4040 Linz, Austria }
\email{norayr.matevosyan@oeaw.ac.at}
\title[Free boundary near contact points]{Behavior of the free boundary near
contact points with the fixed boundary for nonlinear elliptic equations}
\maketitle

\begin{abstract}
The aim of this paper is to study a free boundary problem
for a uniformly elliptic fully non-linear operator. Under certain assumptions we show that
free and fixed boundaries meet tangentially at contact points.
%Draft version \today.
\end{abstract}

\section{Introduction and main results}

In this paper we consider a free boundary problem for a uniformly elliptic
fully non-linear operator $F$ in the following setting:%
\begin{equation}
\left\{
\begin{array}{ll}
F\left( D^{2}u\right) =\chi _{\Omega } & \text{in }B_{1}^{+},\text{ for an
open set }\Omega=\Omega(u) \subset B_{1}^{+} \text{ defined by } \\
u=|\nabla u|=0 & \text{in}\ B_{1}^{+}\setminus \Omega , \\
u=0 & \text{on }\Pi :=\{x_{1}=0\}.%
\end{array}%
\right.  \label{pr}
\end{equation}
where $n \geq 2$ and the PDE holds in the viscosity sense:

%Conditions on $F$
The following conditions are imposed on $F$ throughout the paper:

\begin{enumerate}
\item $F$ is uniformly elliptic with ellipticity constants $\lambda $ and $%
\Lambda $, i.e.,
\begin{equation*}
\lambda \Vert N\Vert \leq F(A+N)-F(A)\leq \Lambda \Vert N\Vert
\end{equation*}%
where $A$ and $N$ are arbitrary $n\times n$ symmetric matrices with $N\geq 0$%
.

\item $F$ is homogeneous of degree one, i.e.,
\begin{equation*}
F(tA)=tF(A),
\end{equation*}%
for all real numbers $t$ and matrices $A$.

\item $F$ is convex.

\item $F$ is $C^{1}$.
\end{enumerate}

\newpage
\begin{definition}
\cite{CC} A continuous function $u$ is a viscosity solution
 of equation
$$
F\left( D^{2}u\right) =\chi _{\Omega }
$$
in $B_{1}^{+}$, when the following condition holds: for any
$x^{0}\in B_{1}^{+}$ and quadratic polinmoials $\varphi$, $\psi $  such that
$u-\phi$ has a local maximum at $x^{0}$ and $u-\psi$ has a local minimum
 at $x^{0}$, the following holds
\begin{eqnarray*}
F\left( D^{2}\varphi \left( x^{0}\right) \right) \geq \chi _{\Omega }\left(x^{0}\right), \\
F\left( D^{2}\psi \left( x^{0}\right) \right) \leq \chi _{\Omega }\left(x^{0}\right).
\end{eqnarray*}
\end{definition}

\noindent Let us denote the free boundary
$$
\{x:u=|\nabla u(x)|=0\}\cap\partial\Omega
$$
by $\Gamma=\Gamma(u)$ and the complement of $\Omega$
$$
\Lambda =\Lambda (u)=B_{1}^{+}\setminus \Omega(u)=\left\{ x\in B_{1}^{+}:u\left( x\right)
=\left\vert \nabla u\left( x\right) \right\vert =0\right\}.
$$

It is well known that viscosity solutions of fully nonlinear uniformly elliptic PDEs have the
usual maximum/minimum principle as well as compactness properties.
Furthermore, they are uniformly $C^{1,\alpha }$ when
$$
|F(D^{2}u)|\leq C,
$$
and $C^{2,\alpha }(B(0,1/2))$ when
$$
F(D^{2}u)=1 \quad \text{in} \ B(0,1).
$$
For the details we refer to \cite{CC}.\newline

In the future we shall use the following notations:
\begin{equation*}
\begin{array}{ll}
\mathbf{R}_{+}^{n} & \qquad \{x\in \mathbf{R}^{n}:x_{1}>0\}, \\
\mathbf{R}_{-}^{n} & \qquad \{x\in \mathbf{R}^{n}:x_{1}<0\}, \\
B(z,r) & \qquad \{x\in \mathbf{R}^{n}:|x-z|<r\}, \\
B^{+}(z,r) & \qquad \{x\in \mathbf{R}_{+}^{n}\cap B(z,r)\}, \\
B^{-}(z,r) & \qquad \{x\in \mathbf{R}_{-}^{n}\cap B(z,r)\}, \\
B_{r},\ B & \qquad B(0,r),\quad B_{1}, \\
B_{r}^{+},\ B^{+} & \qquad B^{+}(0,r),\quad B_{1}^{+}, \\
B_{r}^{-},\ B^{-} & \qquad B^{-}(0,r),\quad B_{1}^{-}, \\
\Pi ,\ \Pi (z,r),\ \Pi _{r} & \qquad \{x\in \mathbf{R}^{n}:x_{1}=0\},\quad
\Pi \cap B(z,r),\quad \Pi (0,r), \\
\Vert \cdot \Vert _{\infty } & \qquad \text{canonical norm}, \\
e_{1},\ldots ,e_{n} & \qquad \text{standard basis in}\ \mathbf{R}^{n}, \\
\nu ,\ e & \qquad \text{arbitrary unit vectors}, \\
D_{\nu },\ D_{\nu e} & \qquad \text{first and second directional derivatives}%
, \\
v_{+},\ v_{-} & \qquad \max (v,0),\ \max (-v,0), \\
\chi _{D} & \qquad \text{the characteristic function of the set}\ D, \\
\partial D & \qquad \text{the boundary of the set}\ D, \\
\Omega=\Omega \left( u\right) & \qquad B^{+}\setminus \left\{ x:u(x)=|\nabla u(x)|=0\right\} , \\
\Omega ^{+}\left( u\right) & \qquad \left\{ x\in \Omega :u\left( x\right)
>0\right\} , \\
\Omega ^{-}\left( u\right) & \qquad \left\{ x\in \Omega :u\left( x\right)
<0\right\} , \\
\Lambda =\Lambda (u) & \qquad \left\{ x\in B_{1}^{+}:u\left( x\right)
=\left\vert \nabla u\left( x\right) \right\vert =0\right\} , \\
\Gamma =\Gamma (u) & \qquad \{x:u=|\nabla u(x)|=0\}\cap \partial \Omega
\text{ the free boundary,} \\
\Gamma ^{\ast }(u) & \qquad \Gamma (u)\cap \Pi \text{ the set of contact
points.}%
\end{array}%
\end{equation*}

We define the density function $V_{r}$ as
\begin{equation*}
V_{r}(z,u)={\frac{\text{vol}(\Omega ^{-}(u)\cap B^{+}(z,r))}{r^{n}}}\,.
\end{equation*}%
Observe that by continuity of $u$, if $V_{r}(0,u)=0$ for all $r$, then $u\geq 0$.

Assume that the set $\{u<0\}$ is small enough near the origin, we shall
prove the
quadratic growth of $|u|$ near the origin (Theorem A). But since in Theorem
A we only assume $V_{r}(z,u)\leq C_{0}$ for $z=0$, we cannot prove that $u$
is $C^{1,1}$ near the origin. However, Theorem~A still allows us to re-scale
$u$ quadratically such that ($(u(rx)/r^{2}$) and remain bounded as $r$ tends to
zero.

In Theorem B we give the classification of positive global solutions
(solutions in the half space $\mathbb{R}_{+}^{n}$). The latter we use in the
proof of Theorem C, which concerns the tangential approach of the free
boundary $\Gamma $ to the fixed boundary $\Pi $. In the special case of
Theorem C, global solutions obtained by the procedure of the blowup (see
below) are nonnegative due to the zero-density assumption on the set $%
\{u<0\} $ at the origin.

The main difficulty in proving this result in the nonlinear case,
without any a priori conditions on the solution, lies in the lack of the
monotonicity lemma, which is heavily used in the case of the Laplacian.

To deal with the complications described above, we use the technique
developed in [LS].

\begin{definition}[Local solution]
We say a continuous function $u$ belongs to the class $P_{r}^{+}(z,M)$ if $u$
satisfies

\begin{enumerate}
\item $F(D^{2}u)=\chi _{\Omega }$ in $B^{+}(z,r)$ in the viscosity sense, for some open set $\Omega$,

\item $u=|\nabla u|=0$ in $B^+(z,r)\setminus \Omega $,

\item $\Vert u\Vert _{\infty ,B^+(z,r)}\leq M$,

\item $u=0$ on $\Pi_1 ,$

\item $z\in \partial \Omega $.
\end{enumerate}
\end{definition}

\begin{definition}[Global solutions]
We say a continuous function $u$ belongs to the class $P_{\infty }^{+}(z,M)$
if $u$ satisfies

\begin{enumerate}
\item $F(D^{2}u)=\chi _{\Omega }$ in $\mathbb{R}_{+}^{n}$ in the viscosity
sense, for some open set $\Omega$,

\item $u=|\nabla u|=0$ in $\mathbb{R}_{+}^{n}\setminus \Omega $,

\item $|u(x)|\leq M(|x|^{2}+1)$,

\item $u=0$ on $\Pi ,$

\item $z\in \partial \Omega $.
\end{enumerate}
\end{definition}

Our first result asserts that local solutions have quadratic growth if the
set $\Omega ^{-}$ is sufficiently small. See \cite{CKS,KS} for a similar type
of result for the Laplacian case. \vspace{3mm}

\noindent \textbf{Theorem A}. \textsl{\ \ There is a universal constant $%
C_{0}=C_{0}(n,\lambda ,\Lambda )$ such that, for all $r<1$ and $u\in
P_{1}^{+}(z,M)$ we have }
\begin{equation*}
\sup_{B^{+}(z,r)}|u|\leq \frac{M}{C_{0}}r^{2}\ \text{for }r<1,
\end{equation*}%
\textsl{\ provided $V_{r}(z,u)\leq C_{0}$, for all $r<1$.}

\begin{remark}
Let $d(x,\partial \Omega )$ be the distance from $x$ to $\partial \Omega .$
From Theorem~A we have that if $V_{r}(y,u)\leq C_{0}$ for all $y\in
B^{+}(0,1/2)\cap \partial \Omega $, and $r<1/2$, then for $u\in P_1(0,M)$
\begin{equation*}
|u(x)|\leq \frac{M}{C_{0}}d^{2}(x,\partial \Omega ).
\end{equation*}%
The interior $C^{2,\alpha }$ estimates (see \cite{CC}) coupled with the
above remark imply that $u\in C^{1,1}\left( B^{+}(0,\frac{1}{4})\right) $
uniformly for the class $P_{1}^{+}(0,M),$ provided the assumption $%
V_{r}(y,u)\leq C_{0}$ holds for all $y\in B^{+}(0,1/2)\cap \partial \Omega $%
, and $r<1/2$. For details we refer to the proof of theorem 1.1 in \cite{CKS}%
.
\end{remark}

\vspace{3mm}

\noindent \textbf{Theorem B}. \textsl{\ \ Let $u\in P_{\infty }^{+}(0,M)$,
\begin{equation*}
u\geq 0, \quad \hbox{and}\quad |\nabla u(0,x_2,\cdots ,x_n)|=0.
\end{equation*}
Then $\Lambda (u)=\Pi$.}

In order to state our last theorems we will need the following definitions.

\begin{definition}
Let $\sigma (r)(\sigma (0^{+})=0)$ be a modulus of continuity. Then we
define $P_{r}^{+}(0,M,\sigma )$ as the subset of all functions $u$ in $%
P_{r}^{+}(0,M)$ with the property
\begin{equation*}
V_{r}(0,u) \leq \sigma (r) \qquad \hbox{and}\qquad |\nabla u(0,x_2,\cdots ,
x_n)| \leq |x|\sigma (|x|).
\end{equation*}
\end{definition}

We believe that both these properties are superfluous and they could be
relaxed. However, at this moment, the lack of techniques such as
monotonicity formulas forces us to restrict ourselves to such cases.

The idea with these assumption is that after scaling the limit functions
will be non-negative, and the gradient will be zero on the fixed boundary.
\vspace{3mm}

% Theorem 3
\noindent \textbf{Theorem C}. \textsl{\ \ Let $M>0$ and $\sigma $ be modulus
of continuity. Then there exists $r_0=r_0(n,M)>0$ and a modulus of
continuity $\sigma_{1}(\sigma_{1}(0^+)=0)$ such that if $u\in
P_{1}^{+}(0,M,\sigma)$, then }
\begin{equation}
\partial\Omega\cap B_{r_0}\subset\{x:x_1\leq\sigma_1(|x|)|x|\}.
\label{tnagentiality}
\end{equation}

The interior case of problem (1.1) (i.e. the problem in the whole ball $B_{1}$%
) has been considered earlier in [LS]. When $F$ is the Laplacian operator,
the problem is considered in [SU], and for the parabolic
operator in [ASU1], [ASU2], [CPS]. See also the pioneering work
of L. A. Caffarelli [C1].
%
%
%
%%%%%% End of Introduction %%%%%%%%%%%%%%%%

\section{Proof of Theorem A}

We will show that, the solution $u$ grows away from the free boundary at
most with a quadratic rate. We follow the main idea given in [LS], which is
to use a homogeneous stretching of the solution by the maximum of $u$ over
the ball $B^{+}(0,r)$. Then we will have a control over the growth of these
functions, and we can consider their limit as $r$ tends to zero.

We define
\begin{equation*}
S_{j}(z,u)=\sup_{B^{+}(z,2^{-j})}|u|.
\end{equation*}
In view of the results of [LS] it will be sufficient to
to prove the following lemma.

\begin{lemma}
There exist a constant $C_{0}$ depending only on $n$, such that for every $%
u\in P_{1}^{+}(z,M)$, $j\in \mathbb{N}$ and $z\in\Gamma(u)\cap B_{1/2}$
\begin{equation}
S_{j+1}(z,u)\leq max\{S_{j}(z,u)2^{-2},C_{0}M2^{-2j}\}
\end{equation}
\label{LforN} provided
\begin{equation*}
V_{2^{-j}}(z,u)\leq C_{0}.
\end{equation*}
The constant $C_{0}$ depends on $n$, $\lambda $, and $\Lambda $.
\end{lemma}

\proof If the conclusion in the lemma fails, then there exist sequences
\begin{equation*}
\{\Omega _{j}\},\{u_{j}\}\subset ~P_{1}^{+}(0,M),\ \{z_j\}\subset\Gamma(u_j)\cap B_{1/2} \ \{k_{j}\}\subset \mathbb{N%
},\ k_{j}\nearrow \infty
\end{equation*}
such that
\begin{equation*}
S_{k_{j+1}}(z_j,u_{j})>max\{2^{-2}S_{j}(z_j,u_{j}),\ Mj2^{-2k_{j}}\}\qquad \forall
j\in \mathbb{N}.
\end{equation*}

Consider the following scaling
\begin{equation*}
\tilde{u}_{j}(x)=\frac{u_{j}(z_j+2^{-k_{j}}x)}{S_{k_{j}+1}(z_j,u_{j})}\qquad \text{%
in $B_{1}^{+}$}.
\end{equation*}%
The following results can be obtained by computation:

\begin{itemize}
\item $\Vert \tilde{u}_{j}\Vert _{\infty ,B}=\frac{S_{k_{j}}(z_j,u_{j})}{%
S_{k_{j}+1}(z_j,u_{j})}\leq 4$ ,

\item $\|\tilde{u}_j\|_{\infty,B_{\frac1{2}}}=1$ ,

\item $\tilde{u}_j(0)=|\nabla\tilde{u}_j(0)|=0$,

\item $V_{1}(\tilde{u}_{j})\leq \frac{1}{j}\to 0$.
\end{itemize}

Also, as in [LS], by ellipticity and degree one homogeneity of $F$ (where $F
$ itself may vary within the bounds of the condition stated earlier)
\begin{equation*}
|F(D^{2}\tilde{u}_{j}(x))|\leq \Lambda \frac{(2^{-k_{j}})^{2}}{%
S_{k_{j}+1}(z_j,u_{j})}\leq \frac{\Lambda S_{k_{j}}(z_j,u_{j})}{jMS_{k_{j}+1}(z_j,u_{j})}%
\leq \frac{4\Lambda }{jM}\to 0\,.
\end{equation*}

Standard elliptic estimates \cite{CC} imply a uniform bound for the $%
C^{1,\alpha }$-norms of $\tilde{u}_{j}$. Therefore a subsequence of $\{%
\tilde{u}_{j}\}$ converges to a function $\tilde{u}_{0}$ satisfying
\begin{eqnarray*}
F(D^{2}\tilde{u}_{0})=0\ \hbox{in }B^{+}(0,1),\ \ \tilde{u}_{0}\geq 0, \\
\tilde{u}_{0}(0)=|\nabla \tilde{u}_{0}(0)|=0\,\ \ \hbox{and}\ \sup_{B_{1/2}}%
\tilde{u}_{0}=1.
\end{eqnarray*}%
The above in particular implies that the nonzero solution $\tilde{u}_{0}$ of
the elliptic equation $F(D^{2}\tilde{u}_{0})=0$ has a local minimum at a
boundary point and its gradient is zero at that point. Using Hopf type lemma
we come to a contradiction.

\section{Nondegeneracy}

We will be concerned with scaling of the type
\begin{equation*}
u_r(x):=\frac{u(rx)}{r^2},
\end{equation*}
and its limit (when it exists)
%The uniform limit of $u_{r_j}$ when $r_j\to0$ is called blow-up of $u$.%
\begin{equation*}
u_0:=\lim_{r_j \to 0} u_{r_j},
\end{equation*}
called blow-up limit. Hence we need to assure that $u_0\neq 0$, i.e. $u$ is
non-degenerate.

\begin{lemma}
If $u\in P_R^+(z,M)$, $x^0\in \overline{ \{u> 0\}}\cap B_{R/2}(z)$ then
%for $r<dist(x_0,\partial B)$
\begin{equation}
\sup_{B^+(x^0,r)}{u}\geq u(x^0)+C_{0}r^2, \qquad \text{for all }\ r<R-|x_0
-z |,  \label{eq-non-deg}
\end{equation}
where $C_{0}$ is a constant depending only on $n$ and $\Lambda$. \label%
{lemma-non-deg}
\end{lemma}

\proof It suffices to consider the case $x^{0}\in \{u>0\} \cap B_{R/2}(z)$
because if \eqref{eq-non-deg} holds for all $x^{0}\in \{u>0\} \cap
B_{R/2}(z) $, then it will be true also for all $x^{0}\in \overline{\{u>0\} }
\cap B_{R/2}(z)$. Set
\begin{equation}
v(x)=u(x)-u(x^{0})-\frac{1}{2n\Lambda }|x-x^{0}|^{2}.  \label{v(x)}
\end{equation}
There exists $x^{1}\in \overline{B^{+}}(x^{0},r)$ such that the following
holds:
\begin{equation}
v(x^{1})=\sup_{B^{+}(x^{0},r)}{v}.  \label{v(x^1)}
\end{equation}%
To prove the lemma, it is enough to prove the following two steps:

\begin{itemize}
\item $v(x^1)\geq0$,

\item $x^1\in\partial B^+(x^0,r)\backslash\Pi(x^0,r)$.
\end{itemize}

\vspace{1mm} The first step simply follows from the fact that
\begin{equation*}
v(x^{1})\geq v(x^{0})=0.
\end{equation*}
To prove the second step assume $x^{1}\in B^{+}(x^{0},r)$. Then from %
\eqref{v(x^1)} we have $|\nabla v|(x^{1})=0$. Thus by \eqref{v(x)}
\begin{equation*}
(\nabla u)(x^{1})=\frac{1}{n\Lambda }(x^{1}-x^{0}).
\end{equation*}%
Now, if $x^{1}\neq x^{0}$, then $(\nabla u)(x^{1})\neq 0$, i.e., $x_{1}\in
\Omega $. We also have
\begin{equation*}
F(D^{2}v)=F\left( D^{2}u-\frac{I}{n\Lambda }\right) \geq F(D^{2}u)-\Lambda
\frac{1}{\Lambda }=0\ \hbox{in }\Omega,
\end{equation*}%
and \eqref{v(x^1)} together with maximum principle gives us that
\begin{equation*}
v(x)\equiv constant=:C \quad \text{in } \overline{\Omega \cap B^{+}(x^{0},r)}%
.
\end{equation*}
In particular, $C=v(x^{0})=0$ so we have
\begin{equation*}
u(x)=u(x^{0})+\frac{1}{2n\Lambda }|x-x^{0}|^{2}
\end{equation*}
and
\begin{equation*}
(\nabla u)(x)=\frac{1}{n\Lambda }(x-x^{0}) \quad \text{in } \overline{\Omega
\cap B^{+}(x^{0},r)}.
\end{equation*}
But if we take $y\in \partial \Omega \cap \overline{B^{+}(x^{0},r)}$ (we may
assume it exists without loss of generality) then we get
\begin{equation*}
|\nabla u(y)|=\frac{1}{n\Lambda }(y-x^{0})\neq 0,
\end{equation*}
which is a contradiction, since $|\nabla u|=0$ on $\partial \Omega $ . So in
this case $x^{1}\in\partial B^{+}(x^{0},r)$.

If $x^1=x^0$, then again $x^1=x^0\in\Omega$ and we have the same
contradiction as above.

Finally, if $x^{1}\in \Pi (x^{0},r)$, then because $u(x^{0})\geq 0$, we get
the following contradiction
\begin{equation*}
0>v(x^{1})\geq v(x^{0})=0,
\end{equation*}
where the second inequality follows from the definition of $x^{1}$. \qed

\section{Proof of Theorems B and C}

\noindent \textbf{Proof of Theorem B}\newline

Under the conditions imposed on $u$, more exactly $\nabla u=0$ on $%
\{x_{1}=0\}$, one can give a proof of Theorem B by continuing the function $u
$ as zero to the lower half space $\mathbb{R}_{-}^{n}$ to obtain a solution
in whole $\mathbb{R}^{n}$ (one can show that in the viscosity sense there is
no mass on $\{x_{1}=0\}$, since there is no jump in the gradient). Then from
the interior result [LS] it follows that the coincidence set is convex,
hence we have a halfspace solution.\\

\noindent For completeness we give a detailed proof based on ideas of [LS].\\

From the convexity of $F$ in $\mathbb{R}_{+}^{n}$ we can conclude that $%
D_{ee}u$ is a supersolution to the linearized problem and hence it has the
minimum principle (see [LS]). We will prove that $u$ is convex using a
contradictory argument. Assume there is a direction $e$ such that
\begin{equation*}
-\infty <\inf_{\Omega (u)}D_{ee}u=-C<0\,.
\end{equation*}%
Then there exists a sequence $\{x^{j}\}$ such that
\begin{equation*}
D_{ee}u(x^{j})\to -C\ \mbox{as}\ j\to +\infty \,.
\end{equation*}

Let us consider the blowup of $u$ with $d_j=\text{dist}(x^j,\partial
\Omega)<+\infty$
\begin{equation*}
u_j(x)=\frac{u(x^j+d_jx)}{d_j^2}\,.
\end{equation*}
We remark that by the assumption $|\nabla u (0,x_2\cdots ,x_n)|=0$ and
Theorem A (since $u\geq 0$) we have $u_r$ is uniformly bounded. Using
compactness argument we get
\begin{equation*}
u_j\to u_0 \ \ \hbox{   in } C^{2,\alpha}(B^+_{1/2})
\end{equation*}
which implies
\begin{equation*}
D_{ee}u_0(0)= \lim_j D_{ee}u_j(0) = \lim_j D_{ee}u(x^j) =-C
\end{equation*}
and
\begin{equation*}
D_{ee} u_j(x) = D_{ee} u (x^j + d_jx) \geq -C\,.
\end{equation*}
Thus in $B^+_{1/2}$ we have
\begin{equation*}
D_{ee}u_0(x)\geq -C.
\end{equation*}
By maximum principle $D_{ee}u_0\equiv -C$ in $\Omega^{\prime}$, the
connected component of $\Omega (u_0)$ containing the origin. Following the
steps in [LS], we rotate the coordinate system such that $e$ coincides with $%
e_1$.
%and obtain that $D_1 u_0 $ is non-increasing  in the $x_1$-direction,
Next we integrate $D_{11} u_0$ and use non-negativity of %$$
%u_0 = -Cx^2_1/2 + g(x_1^{\prime})x_1+ f(x_1^{\prime}),
%$$
%where $x_1^{\prime}=(x_2,\cdots , x_n)$. Since
$u_0 \geq 0$ to obtain
%we must have $x^2_1 \leq (g(x_1^{\prime})x_1 + f(x_1^{\prime}))/C$ for $x \in\Omega^{\prime}$.
$|x_1| \leq G(x_1^{\prime})$ for some function $G$ and all $x \in
\Omega^{\prime}$. Now, for fixed $x$ let us consider $x^m:= (x_1+m, x_2
\cdots , x_n)$. There exists $m$ depending on $x_1^{\prime}$ such that
\begin{equation*}
u_0(x^m)=|\nabla u_0(x^m)|=0.
\end{equation*}
Combining the facts that for large $m$ we have $D_1 u_0(x^m)=0$, and $D_1 u_0
$ is non-increasing, we get
\begin{equation*}
D_1 u_0 \leq 0 \text{ in } \Omega^{\prime}.
\end{equation*}
The latter gives a contradiction to the non-degeneracy, Lemma 3.1.

Now, as we have proved the convexity of $u$ and hence of the set $\Lambda
=\{u=0\}$, we see clearly that $\Lambda (u) =\Pi$. Indeed, since $\Pi
\subset \Lambda$, if there is a point $x^0\in \Lambda \cap \mathbf{R}^n_+$,
then by convexity
\begin{equation*}
\{x:\quad 0\leq x_1\leq x_1^0\}\subset \Lambda,
\end{equation*}
implying that the origin is not a free boundary. This is a contradiction. %
\qed

\noindent \textbf{Proof of Theorem C}\newline
It is enough to check that for every given $\epsilon$ there exists $%
\rho=\rho_\epsilon$ such that for all $x^0\in\partial\Omega\cap
B^+_{\rho_\epsilon}$
\begin{equation}
x^0\in B^+_{\rho_\epsilon}\backslash K_\epsilon,  \label{tang}
\end{equation}
where
\begin{equation*}
K_\epsilon=\{x:x_1>\epsilon(x_2^2+\ldots+x_n^2)^{1/2}\}.
\end{equation*}
Then we may choose $r_0=\rho_{\{\epsilon=1\}}$ and $\sigma$ given by the
inverse of $\epsilon\to\rho_\epsilon$. The proof is based on a contradictory
argument. If \eqref{tang} fails, than there exists a sequence
\begin{equation*}
u_j\in P_1^+(0,M,\sigma),\ x^j\in\partial\Omega(u_j)\cap B^+_{\rho_j}
\end{equation*}
such that $\rho_j\to0$ and $x^j\in B^+_{\rho_j}\cap\overline{K_\epsilon}$.
Now for every scaled function
\begin{equation*}
\tilde{u}_j(x)=u_j(x|x^j|)/|x^j|^2
\end{equation*}
we have a point $\tilde{x} ^j\in \partial B_1^+\cap \partial\Omega(\tilde{u}%
_j)\cap K_\epsilon$. There exists converging subsequences of $\tilde{u}_j\to
u_0$ and $\tilde{x}^j\to x_0$ such that $x_0\in\overline{K}%
_\epsilon\cap\partial B_1$, with $u(x^0)=0$. It follows by the assumption
\begin{equation*}
|\nabla u_j(0,x_2,\cdots ,x_n)|\leq |x|\sigma (|x|),
\end{equation*}
that
\begin{equation*}
|\nabla \tilde u_j(0,x_2,\cdots ,x_n)|\leq |x|\sigma (|x||x^j|)\quad \to
\quad 0.
\end{equation*}
In particular (by non-degeneracy lemma) $u_0$ is a nonzero global solution,
satisfying the assumptions of Theorem B, and hence $\Lambda =\Pi$,
contradicting $x^0 \in \partial \Omega (u_0)$. \qed

\end{document}